\documentclass[draft]{lematema2}
\usepackage{latexsym}
\usepackage{amsfonts}
\usepackage{amssymb}
\usepackage{amssymb}
\usepackage{amscd}

\theoremstyle{plain}
\newtheorem{theorem}{Theorem}[section]
\newtheorem*{theorem*}{Theorem}

\newtheorem{corollary}[theorem]{Corollary}
\newtheorem*{corollary*}{Corollary}

\newtheorem{fact}[theorem]{Fact}

\newtheorem{properties}[theorem]{Properties}

\newtheorem*{proposition*}{Proposition}

\theoremstyle{definition}

\newtheorem{definition}[theorem]{Definition}
\newtheorem{definitions}[theorem]{Definitions}

\newtheorem*{example*}{Example}

\theoremstyle{remark}

\DeclareMathAlphabet{\mathcal}{OMS}{cmsy}{m}{n}

\def\PV{{\mathbb{P}(V)}}
\def\Pe{{\mathbb{P}^{2n+1}}}
\def\Ou{\mathcal{O}}

\DeclareMathOperator{\im}{im}
\DeclareMathOperator{\rk}{rk}

\title{On the non-existence of orthogonal instanton bundles on $\Pe$}

\titlemark{On the non-existence of orthogonal instanton bundles on $\Pe$}

\titlenote{}
\titlenote{}


\author{\L ucja Farnik}
\address{%
Instytut Matematyki\\
Uniwersytet Jagiello\'nski\\
ul. \L ojasiewicza 6\\
30-348 Krak\'{o}w (Poland)\\
\email{lucja.farnik@im.uj.edu.pl}
}
\author{Davide Frapporti}
\address{%
Dipartimento di Matematica\\
Universit\`a degli Studi di Trento\\
via Sommarive, 14\\
38123 Povo TN (Italy)\\
\email{frapporti@science.unitn.it}
}
\author{Simone Marchesi}
\address{%
Dipartimento di Matematica ``Federigo Enriques''\\
Universit\`a degli Studi di Milano\\
via Cesare Saldini, 50\\
20133 Milano (Italy)\\
\quad\\
Departamento de \'Algebra\\
Facultad de Ciencias Matem\'aticas\\
Universidad Complutense de Madrid\\
Plaza de las Ciencias, 3\\
28040 Madrid (Spain)\\
\email{simone.marchesi@unimi.it\\
smarches@mat.ucm.es}
}
\date{2009/12/11}

\begin{document}

\maketitle

\thispagestyle{empty}
%$\math$

\begin{abstract} In this paper we prove that there do not exist orthogonal instanton bundles on $\Pe$. In order to demonstrate this fact, we propose a new way of representing the invariant, introduced by L.~Costa and G.~Ottaviani, related to a rank $2n$ instanton bundle on $\mathbb{P}^{2n+1}$.
\end{abstract}

\section{Introduction}
Instanton bundles were introduced in \cite{ADHM} on $\mathbb{P}^3$ and in \cite{okon-spin} on $\Pe$. Recently many authors have been dealing with them. The geometry of special instanton bundles is discussed in \cite{spin-traut}.
Authors of \cite{anc-ott} study the stability of instanton bundles. They prove that all special symplectic instanton bundles on $\Pe$ are stable and all instanton bundles on $\mathbb{P}^5$ are stable. Some computations involving instanton bundles can be done using computer programs, e.g. Macaulay, as it is presented in \cite{AABOP}.
The $SL(2)$-action on the moduli space of special instanton bundles on $\mathbb P^3$ is described in \cite{costa-ott2}.
In \cite{costa-ott} it is shown that the moduli space of symplectic instanton bundles on $\Pe$ is affine, by introducing the invariant that we use in this paper.

As we see, instanton bundles were discussed in many papers. 
On $\mathbb{P}^3$ they have rank $2$ and they are symplectic.
On $\Pe$, for $n\ge 2$, theoretically they can be divided in three groups: symplectic, orthogonal and neither symplectic nor orthogonal bundles. 
Lots of examples of symplectic bundles are known as well as of bundles from the third group. 
However, it was not known whether there exist orthogonal instanton bundles. We prove, see Theorem \ref{mainthm}, that it's impossible for such bundles to exist. The question was posed by Trautmann in \cite{traut}.

The paper is organised as follows. In the Preliminaries we give definitions of instanton bundle, orthogonal instanton bundle and symplectic instanton bundle. We state equivalent conditions for a bundle to be respectively orthogonal or symplectic. In the main part we prove the fact that there do not exist orthogonal instanton bundles on $\Pe$ for an arbitrary $n$ and arbitrary $k=c_2(E)$. 

\section{Preliminaries}

\begin{definition}
Let $X$ be a smooth projective variety. A \textit{monad} on $X$ is a complex of vector bundles:
$$0\longrightarrow F \stackrel{\beta}{\longrightarrow}G \stackrel{	\alpha}{\longrightarrow}H \longrightarrow 0$$
which is exact at $F$ and at $H$.
\end{definition}

Let $V$ be a vector space of dimension $2n+2$ and $\Pe=\PV$. Let $U$ and $I$ be vector spaces of dimension $k$ and let $W$ be a vector space of dimension $2n+2k$. We consider the following monad on $\Pe$:
\begin{equation}
U \otimes \Ou_{\Pe}(-1) \stackrel{\beta}\longrightarrow W\otimes \Ou_{\Pe} \stackrel{\alpha}\longrightarrow I \otimes \Ou_{\Pe}(1)  .
\label{monad}
\end{equation}

Note that $\alpha$ and $\beta$ can be represented by two matrices $A$ and $B^t$ respectively, where $A$ and $B$ are $k\times(2n+2k)$-matrices with linear entries of maximal rank for every point in $\Pe$ and $A \cdot B^t=0$.
The condition for $A$ and $B$ to be of maximal rank for every point in $\Pe$ is equivalent to require that the map $\alpha$ is surjective and the map $\beta$ is injective.

\

From now on, with a little abuse of notation,
we will identify a linear map $\alpha$ with the matrix $A$ associated to it, and we will simply write $A$.

\begin{definition}
We say that $E$ is an \textit{instanton bundle} (with quantum number $k$) if $E$ is the cohomology of a monad
as in (\ref{monad}), i.e.
$$E=\ker A/ \im B^t.$$
\end{definition} 

\begin{properties}\quad 
\begin{itemize}
\item[--] For an instanton bundle $E$ we have that
$\rk E=2n$.
\item[--] The Chern polynomial of $E$ is $c_t=(1-t^2)^{-k}=1+kt^2+ \dbinom{k+1}{2}t^4+\ldots$\\
In particular, $c_1(E)=0$ and $c_2(E)=k$. 
\item[--] $E(q)$ has natural cohomology in the range $-2n-1 \leq q \leq 0$, i.e.\\ $h^i(E(q))\neq 0$ for at most one $i=i(q)$.
\end{itemize}
\end{properties}

\begin{definitions}\label{def:ort}
We say that an instanton bundle $E$ is \textit{symplectic} if there exists an isomorphism $\alpha: E\longrightarrow E^*$ such that $\alpha=-\alpha^*$.

We say that an instanton bundle $E$ is \textit{orthogonal} if there exists an isomorphism $\alpha: E\longrightarrow E^*$ such that $\alpha=\alpha^*$.
\end{definitions}

Having a symplectic instanton bundle is equivalent to the existence of a non-degenerate, skew-symmetric matrix $J$ such that $B=A\cdot J$. After linear change of coordinates we may assume that $J$ is of the form

$$J= \left(\begin{matrix} 0 & I\\ -I & 0
\end{matrix}\right)$$
and $A\cdot J\cdot A^t=0$.

Having an orthogonal instanton bundle is equivalent to the existence of a non-degenerate, symmetric matrix $J$ such that $B=A\cdot J$. We may assume that $J=I$
and $A\cdot J \cdot A^t=A\cdot A^t=0$.

\section{On the non-existence}
Recalling the definition we gave before, having an orthogonal instanton bundle means having an application
\begin{equation}\label{def:1}
W \otimes \mathcal{O}_{\mathbb{P}^{2n+1}} \stackrel{A}{\longrightarrow} I \otimes \mathcal{O}_{\mathbb{P}^{2n+1}}(1)
\end{equation}
with the additional properties that $A\cdot A^t = 0$, and $A$ has maximal rank for every point in $\Pe$.

%We also recall that $\mbox{dim}W = 2n+2k$, $\mbox{dim}I = k$ and $\mbox{dim}V = 2n+2$, where $\mathbb{P}\left(V\right) = \mathbb{P}^{2n+1}$.

Let us choose a basis for each vector space we are considering. In particular we have $i_1$, $\ldots$, $i_k$ as a basis of $I$, $w_1$, $\ldots$, $w_{2n+2k}$ as a basis of $W$ and $v_0$, $\ldots$, $v_{2n+1}$ as a basis of $V$.

From the map defined in (\ref{def:1}) we can induce a second map
\begin{equation}
W \stackrel{M}{\longrightarrow} I \otimes V
\end{equation}
over the global sections where $M$ is a $\left(k\left(2n+2\right)\right) \times \left(2n+2k\right)$-matrix.

We want to construct such matrix in a particular way, i.e.
$$
M: = \left( \begin{array}{c}
M_1 \\ \hline \vdots \\ \hline M_k
\end{array} \right)
$$
where, fixed $j$, the block $M_j$ is a $\left(2n+2\right) \times \left(2n+2k\right)$-matrix ``associated'' to the element $i_j$ of the basis of $I$. In other words the block $M_j$ represents the following correspondence:
$$
M_j : w \mapsto i_j \otimes M_j(w) \:\:\: \mbox{with} \:\: i_j \:\: \mbox{fixed}.
$$

Using the matrix $M$ it is possible to reconstruct the quadratic conditions given by $A \cdot A^t = 0$. Indeed, let us observe that the matrix $A$ can be described by the following matrix
$$
A = \left( \begin{array}{c}
\left[x_0, \ldots, x_{2n+1}\right] \cdot M_1 \\ \hline \vdots \\ \hline  \left[x_0, \ldots, x_{2n+1}\right] \cdot M_k
\end{array} \right),
$$
where $x_i$ are the coordinates in $V$ with respect to the basis $\{v_0$, $\ldots$, $v_{2n+1}\}$.

The most natural idea is to consider the product $M \cdot M^t$, that %considering the description we gave of the matrix $M$, 
is of the form
$$
M \cdot M^t = \left(
\begin{array}{cccc}
M_1\cdot M_1^t & M_1\cdot M_2^t & \cdots &M_1\cdot M_k^t \\
M_2\cdot M_1^t & M_2\cdot M_2^t & \cdots & M_2\cdot M_k^t \\
\vdots & & &\vdots\\
M_k\cdot M_1^t & M_k\cdot M_2^t & \cdots & M_k\cdot M_k^t \\
\end{array} \right),
$$
where each block $M_\alpha \cdot M_\beta^t$, for $\alpha$, $\beta=1$, $\ldots$, $k$, is a square matrix representing a quadratic form in the variables $x_0$, $\ldots$, $x_{2n+1}$. Because of this considerations, we can state the following equivalence of quadratic conditions
\begin{equation}
A\cdot A^t = 0 \:\:\: \Longleftrightarrow \:\:\: M_\alpha \cdot M_\beta^t + M_\beta \cdot M_\alpha^t = 0, \:\: \mbox{for all} \:\: \alpha,\beta = 1,\ldots,k.
\label{eq:1}
\end{equation}

Let us observe, as in \cite{costa-ott}, that the vector spaces $W \otimes \mathcal{S}^n I$ and $\mathcal{S}^{n+1} I \otimes V$ have the same dimension $\left(2n + 2k \right)\binom{k+n-1}{n} = \left(2n+2\right)\binom{k+n}{n+1}$, where we denote by $\mathcal{S}^n I$ the $n$-th symmetric power of $I$.

We can induce from
$$
W \stackrel{M}{\longrightarrow} I \otimes V
$$
the morphisms
$$
M \otimes id_{\mathcal{S}^{n}I} : W \otimes \mathcal{S}^n I \longrightarrow V \otimes I \otimes \mathcal{S}^n I,
$$
$$
id_V \otimes \pi : V \otimes I \otimes \mathcal{S}^n I \longrightarrow V \otimes \mathcal{S}^{n+1} I,
$$
with $\pi$ denoting the natural projection, and we consider the composition
\begin{equation}
Q = (id_V \otimes \pi) \circ (M \otimes id_{\mathcal{S}^{n}I}).
\end{equation}

We recall a result stated in \cite{costa-ott}, that, in our situation, is the following
\begin{theorem}\label{teo:1}
If we have an instanton bundle given by $A$, then $\det Q \neq 0$, i.e. $Q$ is non-degenerate. $\det Q$ is $SL(W)\times SL(I)\times SL(V)$-invariant.
\end{theorem}

As we did for $M$, we construct the matrix associated to $Q$. First of all we fix the lexicographic order induced by $i_1 >  i_2 > \ldots > i_k$ on the elements of the basis of $\mathcal{S}^{n}I$ and on the elements of the basis of $\mathcal{S}^{n+1}I$.
We label with $\zeta_1,\ldots,\zeta_s$ the elements of $\mathcal{S}^{n}I$ and by $\eta_1,\ldots, \eta_r$ the elements of $\mathcal{S}^{n+1}I$,
where  $s=\binom{k+n-1}{n}$ and $r=\binom{k+n}{n+1}$.
We decompose $Q$ in $r\cdot s$ blocks of dimension $(2n+2) \times (2n+2k)$.
The block $Q_{ij}$,
is the matrix that represents the linear map
$$w\otimes \zeta_j\longmapsto Q_{ij}(w) \otimes \eta_i,$$
for $i\in\{1,\ldots, r\}$ and $j\in\{1,\ldots, s\}$.

We will see later, with more details, that each block $Q_{ij}$ is either one of the $M_j$ blocks or the matrix of all zeros.

%As before, each block $M_j$, for a fixed $j$, represents the application
%$$
%M_j \otimes id_{\mathcal{S}^{n}I} : W \otimes \mathcal{S}^n I \longrightarrow id_V \otimes \pi\left(V \otimes i_j \otimes \mathcal{S}^n I\right).
%$$
\begin{fact}\label{fact:1}
Considering the matrix $Q$ we have constructed, there always exists a syzygy 
 $V^*\otimes \mathcal{S}^{n-1}I$ with %$\mathfrak{S}$
$$
V^*\otimes \mathcal{S}^{n-1}I \stackrel{S}{\longrightarrow} W \otimes \mathcal{S}^n I \stackrel{Q}{\longrightarrow} V \otimes \mathcal{S}^{n+1} I
$$
such that $Q\cdot S = 0$ and $S\neq 0$.
\end{fact}
In particular, we are going to demonstrate that the map $S$, represented by a $\left(\left(2n + 2k \right)\binom{k+n-1}{n}\right) \times \left(2n + 2 \right)$-matrix of the form
\begin{equation}
S = \left( \begin{array}{c}
M_1^t \\ \hline \vdots \\ \hline M_k^t \\ \hline 0 \\ \hline \vdots \\ \hline 0
\end{array} \right)
\label{syzygy}
\end{equation}
will accomplish the property we are looking for.

We are going to abuse some language now, because remembering that $Q$ is divided by blocks, from now on we will refer to a column of blocks of $Q$ as a column of $Q$ and, the same way, to a row of blocks of $Q$ as a row of $Q$.\\
We observe that each column of $Q$ is related to an element of the basis of $\mathcal{S}^n I$ and each row is related to an element of the basis $\mathcal{S}^{n+1} I$. We advice to take a look at the next example to have a better understanding of all the constructions we made.

In order to prove that the map $S$, defined in (\ref{syzygy}), is such that $Q\cdot S=0$, we have only to look at the first $k$ columns of $Q$.
This means that we just need to check all possible passages from the basis elements $\left\{i^n_1,i_1^{n-1}i_2,\ldots,i_1^{n-1}i_k\right\} = \left\{i_1^{n-1}i_\alpha \mid \alpha=1,\ldots,k\right\}=\mathcal{I}\subset \mathcal{S}^n I$ to all the element of the basis of $\mathcal{S}^{n+1} I$ of the form $\left\{i_1^{n-1}i_\alpha i_\beta \mid \alpha,\beta=1,\ldots,k\right\}$. Let us observe that we are only considering all the elements of the two basis that contain the monomial $i_1^{n-1}$.
%\newpage
We now have two options:
\begin{itemize}
\item We can get every element of $\mathcal{S}^{n+1} I$ of the type $i_1^{n-1}i_\alpha^2$, for $\alpha=1$, $\ldots$, $k$, from an element
of $\mathcal{I}$, only multiplying  $i_1^{n-1}i_\alpha$ by $i_\alpha$.\\
This means that the only non-vanishing block of the row of $Q$ corresponding to the element $i_1^{n-1}i_\alpha^2$ will be a block $M_\alpha$ in the $\alpha$-th column.
%This means that the row of $Q$ which corresponds to the element $i_1^{n-1}i_\alpha^2$, for fixed $\alpha$ from 1 to $k$,
% will have an $M_\alpha$ block in the $\alpha$-th column.
 Multiplying this row by $S$ we will obtain $M_\alpha \cdot M_\alpha^t$ which is zero for each $\alpha$, because of the quadratic relations.
 
\item We can get the elements of $\mathcal{S}^{n+1} I$ of the type $i_1^{n-1}i_\alpha i_\beta$, $\alpha \neq \beta$ in only two different ways: multiplying $i_1^{n-1}i_\alpha$ by $i_\beta$ or multiplying $i_1^{n-1}i_\beta$ by $i_\alpha$.\\
This means that the only non-vanishing blocks of the row of $Q$ corresponding to the element $i_1^{n-1}i_\alpha i_\beta $ will be a block $M_\alpha$ in the $\beta$-th column and a block $M_\beta$ in the $\alpha$-th column.
%This means that the row of $Q$ which corresponds to the element $i_1^{n-1}i_\alpha i_\beta$ , for fixed $\alpha,\beta$ from 1 to $k$,
 %will have an $M_\alpha$ block in the $\beta$-th column and an $M_\beta$ block in the $\alpha$-th %column.
 Multiplying this row by $S$ we will obtain $M_\alpha \cdot M_\beta^t + M_\beta \cdot M_\alpha^t $ which is zero for each $\alpha$, $\beta$, because of the quadratic relations.
\end{itemize}
%Let us observe that the remaining part of the matrix $Q$, from the $k+1$-th column till the last one, it's %also constructed by the blocks $M_{\alpha}$, for $\alpha$ from 1 to $k$, and the blocks of zeros. %However, we do not care about such columns because the matrix $S$ we have chosen is zero from the %$k+1$-th row of blocks till the end. This means that the contribution of this last part of $Q$ in the product %$Q \cdot S$ will always be null.\\
Because of everything we said, the Fact \ref{fact:1} is proved.

From Fact \ref{fact:1} we have that the kernel of $Q$ is not zero, so $\det Q =0$. This gives us a contradiction with Theorem \ref{teo:1}.

Due to all our previous work, we can state
\begin{theorem}\label{mainthm}
There do not exist orthogonal instanton bundles $E$ on $\mathbb{P}^{2n+1}$, for each $k$, $n$ $\in$ $\mathbb{N}$.
\end{theorem}

\paragraph{\textbf{Example.}}
We want to write the following simple case in order to give to the reader a better understanding of the previous proof.

%The following lemma gives us a first bound to choose $k$ and $n$ in order to get an orthogonal bundle.
%\begin{lemma}
%Let $E$ be an orthogonal instanton bundle on $\mathbb{P}^{2n+1}$ of quantum number $k$.
%Then $k\geq n+2$. 
%\end{lemma}
%
%\begin{proof}
%The first row of $A=(a_{ij})$ gives an embedding $\varphi: \mathbb{P}^{2n+1}\hookrightarrow\mathbb{P}^{2n+2k-1}$,
%since $A(x)$ has maximal rank for all $x\in \mathbb{P}^{2n+1}$.
%From $A\cdot A^t=0$ we get that  $\varphi(x)\in Q$ for all $x\in \mathbb{P}^{2n+1}$, 
%where $Q$ is the smooth quadric of dimension $2n+2k-2$ in $\mathbb{P}^{2n+2k-1}$ given by $\sum_{j=1}^{2n+2k}a_{1j}^2(x)=0$.
%
%By \cite[Theorem 22.13]{harris-alg} we know that
%on a smooth quadric of dimension $p$, there can be embedded a linear space of 
%dimension at most $\dfrac{p}{2}$. Hence
%$$ 2n+1 \leq \frac{1}{2}(2n+2k-2)$$
%which gives that $k \geq n+2$.
%\end{proof}

It is known that for $n=1$, see \cite{okon-spin}, there do not exist instanton bundles, so taking a simple case, we consider an instanton bundle $E$ with $n=2$  and $k=4$.
Let $i_1$, $\ldots$, $i_4$ be a basis of $I$, $v_0$, $\ldots$, $v_5$ be a basis of $V$ and $w_1$, $\ldots$, $w_{12}$ be a basis of $W$.
 In this particular case we have a map
$$
W \stackrel{M}{\longrightarrow} I \otimes V
$$
with $M$ a $24\times 12$-matrix. Remembering what we said before, $M$ will have the form
$$
M = \left(
\begin{array}{c}
M_1 \\ \hline M_2 \\ \hline M_3 \\ \hline M_4
\end{array}
\right).
$$
%where, for every fixed $j$ from 1 to 4, the block $M_j$ is a $6\times 12$-matrix associated to the application
%$$
%M_j : W \longrightarrow i_j \otimes V.
%$$
%
%From the matrix $M$ we can reconstruct the quadratic conditions $A \cdot A^t = 0$, in fact they are equivalent to the following condition
%\begin{equation}\label{eq:1}
%M_\alpha \cdot M_\beta^t + M_\beta \cdot M_\alpha^t = 0, \:\:\: \alpha,\beta=1,\ldots,4.
%\end{equation}
In this case, % in order to construct the application $Q$, 
we consider the map
$$
W \otimes \mathcal{S}^2 I \stackrel{Q}{\longrightarrow} \mathcal{S}^3 I \otimes V
$$
induced by $M$.
Here $Q$ is represented by a square matrix of order 120. Basing on what we have said, $Q$ will have the following form
$$
Q =\left(
\begin{tabular}{|c|c|c|c|c|c|c|c|c|c||c}
\hline $M_1$&   &  &   &    & & & &    & & $i_1^3$ \\ \hline 
 $M_2$& $M_1$ &   &   &    & & & &    & &    $i_1^2i_2$ \\ \hline
 $M_3$&    &$M_1$ &   &    & & & &    &   & $i_1^2i_3$ \\ \hline
 $M_4$&    &   &$M_1$ &    & & & &    &  &   $i_1^2i_4$ \\ \hline
      &$M_2$  &   &   &$M_1$  & & & &    & & $i_1i_2^2$ \\  \hline
      &$ M_3$ &$M_2$ &   &    &$M_1$ &   &  &    &  & $i_1i_2i_3$\\  \hline
   		&$M_4 $ &   &$M_2$ &    &   &$M_1$ &  &    &  & $i_1i_2i_4$\\  \hline
      &    &$M_3$ &   &    &   &   &$M_1$&    & & $i_1i_3^2$  \\  \hline 
      &    &$M_4$ &$M_3$ &    &   &   &  &$M_1$  & &  $i_1i_3i_4$\\  \hline 
      &    &   &$M_4$ &    &   &   &  &    &$M_1$ &  $i_1i_4^2$ \\ \hline 
      &    &   &   &$M_2$  & & & &    &  & $ i_2^3$\\  \hline
      &    &   &   &$M_3 $ &$M_2$ &   &  &    & & $i_2^2i_3$ \\  \hline
      &    &   &   &$M_4$  &   &$M_2$ &  &    &  & $i_2^2i_4$ \\  \hline
      &    &   &   &    &$M_3$ &   &$M_2$&    &  & $i_2 i_3^2$ \\  \hline
      &    &   &   &    &$M_4$ &$M_3$ &  &$M_2$  &   & $i_2i_3i_4$\\  \hline 
      &    &   &   &    &   &$M_4$ &  &    &$M_2$  & $i_2i_4^2$ \\  \hline
      &    &   &   &    &   &   &$M_3$&    &  & $i_3^3$ \\  \hline
      &    &   &   &    &   &   &$M_4$&$M_3$  &  & $i_3^2 i_4$ \\  \hline
      &    &   &   &    &   &   &  &$M_4$  &$M_3$  & $i_3i_4^2$ \\  \hline
      &    &   &   &    &   &   &  &    &$M_4$& $i_4^3$ \\  \hline\hline
 $i_1^2$ & $i_1i_2$  & $i_1i_3$   & $i_1i_4$   & $i_2^2$ & $i_2i_3$ & $i_2i_4$ & $i_3^2$ & $i_3 i_4$ & $i_4^2$  &
\end{tabular}
\right)
$$
where the empty entries represent $6\times 12$-matrices of zeros. We observe that in the bottom and right part of the matrix we can see, respectively, the elements of the basis of $\mathcal{S}^2 I$ related to the columns of $Q$ and the elements of the basis of $\mathcal{S}^3I$ related to the rows of $Q$.

Taking in mind the quadratic conditions expressed in (\ref{eq:1}), it is immediate to check that
$$
Q \cdot S = Q \cdot \left(
\begin{array}{c}
M_1^t\\
\hline
M_2^t\\
\hline
M_3^t\\
\hline
M_4^t\\
\hline
0\\ \hline
\vdots\\
\vdots\\
 \hline
0\\
\end{array}
\right)
=0,$$
with $S \neq  0$. We can conclude that for the case $k=4$ and $n=2$ there do not exist orthogonal instanton bundles.

\begin{corollary}
Any self-dual instanton bundle $E$ is symplectic.
\end{corollary}
\begin{proof} The bundle $E$ is simple by \cite{anc-ott}. Then 
$$1=h^0(E\otimes E^{*})=h^0(S^2 E)+ h^0(\wedge^2 E).$$  
Theorem \ref{mainthm} implies that $h^0(S^2E)=0$ so we have $h^0(\wedge^2E)=1$ and $E$ is symplectic.
\end{proof}

\section{Acknowledgements}
\quad 
The authors would like to express thanks for the wonderful environment created at P.R.A.G.MAT.I.C. to
the speakers: Rosa Maria Mir\'o-Roig, Giorgio Ottaviani,  Laura Costa and  Daniele Faenzi; and the local organizers:
Alfio Ragusa, Giuseppe Zappal\'a, Rosario Strano, Renato Maggioni and Salvatore Giuffrida, and they also acknowledge the support from the University of Catania during the P.R.A.G.MAT.I.C. 2009 (Catania, Italy, September 2009). 

The authors heartily thank  Giorgio Ottaviani and Daniele Faenzi for many helpful comments and suggestions.
%\newpage
\nocite{*}
\bibhere

\end{document}